
\input vanilla.sty
\magnification 1200
\overfullrule = 0pt
\baselineskip 18pt
\input definiti.tex
\input mathchar.tex
\scriptscriptfont0=\sevenrm
\scriptscriptfont1=\seveni
\scriptscriptfont2 =\sevensy
\scriptscriptfont\bffam=\sevenbf
\define\pmf{\par\medpagebreak\flushpar}

\define\vp{\varphi}
\define\pbf{\par\bigpagebreak\flushpar}

\define\cald{{\cal D}}
\pmf
\title
ALL REGULATORS OF FLAT BUNDLES ARE TORSION
\endtitle
\pmf
\centerline{15 JUNE 1993}
\pmf
\centerline{SECOND REVISION MAY, 1994}
\pmf
\author
Alexander Reznikov
\endauthor
{\it In the Memory of Georgi\u i Isaakovich Katz, the Teacher and the Man.}

\heading{1. Introduction}
\endheading

In our previous paper [22] we proved the Bloch conjecture on
rationality
of secondary characteristic classes of flat rank two vector bundles over a
compact K\"ahler manifold, and announced the full higher dimensional
generalization, which we now state.

\proclaim{1.1 Theorem}(\bf {the Bloch conjecture})  Let $X$ be a smooth complex
projective variety.  Let
$\rho: \pi_1 (X) \to SL_n (\bbc), n \ge ~2$, be a representation of
the fundamental
group and let $E_\rho$ be the corresponding holomorphic flat vector bundle over
$X$.  For $i \ge 2$ let
$c_i \in H^{2i}_\cald (X,\bbz(i))$ be the Chern class in the Deligne cohomology
group of $X$ [10],[12].Then $c_i$ is a torsion class.
\endproclaim

I have some reason to believe that the statement may still be true for
any bundle with a holomorphic
connection. (c.f. [13]).

\pmf
{\bf 1.2} \ The theorem 1.1  contrasts the rank one case, where the
secondary class of the flat bundle, associated to a representation
$\rho : \pi_1 (X) \to GL_1 (\bbc) = \bbc^\ast$ is essentially the
representation itself, so it may be completely arbitrary.  But this
freedom becomes limited when one looks at the
line bundles over curves, corresponding via the Deligne-Ramakrishnan
construction
to elements of $K_2 (X)$.  The map $r : K_2 (X) \to H^1
(X, \bbc^\ast)$ which appears
here, is called the Bloch-Beilinson regulator map.  When $X$ is defined
over a
number field, $k$, the image of $K_2 (X_k)$ under $r$ shows
remarkable properties which were resumed in the
celebrated Bloch-Beilinson conjecture [1], [20].
\pmf
{\bf 1.3} \ The way in which we succeeded to prove the Bloch conjecture for
rank two bundles
in [22]
was the following.  Due to results, proved in succession
by Bloch [2], Soul\'e [26],
and Gillet-Soul\'e [15], one knows that the class
$ c_i$ coincides with the image of the secondary characteristic class
$Ch_{2i-1}(\rho)$ of the flat bundle
$E(\rho)$, lying in $H^{2i-1}(X,\bbc/\bbz)$ and corresponding two the (torsion)
Chern class $c_i(E_\rho)$ in $H^{2i}(X,\bbz)$ under the natural map
$H^{2i-1}(X,\bbc/\bbz)\to H^{2i-1}(X,\bbc/\bbz(i))\to H^{2i}_\cald
(X,\bbz(i))$.
Next, the representation
$\rho$ defines (up to homotopy) a map $\hat \rho: X \to B S
L^\delta_n (\bbc)$, and one knows from Beilinson, Soul\'e, and
Dupont-Sah, that
there exists a universal class $\hat c_{2i-1} \in H^{2i-1} (B S L^\delta_n
(\bbc),
\bbc / \bbz)$, such that $\hat \rho^\ast \hat c_{2i-1} =  Ch_{2i-1}(\rho)$.  An
application of the Vinberg lemma shows that $\hat \rho$ decomposes
through $BSL^\delta_n (\Cal O_S)$,
where $\Cal O_S \subset F$ is a localized ring of integers in a
number field $F$.  With this in mind, the Theorem 1.1 follows from the
following one.
\proclaim{ Main Theorem}  Let $X$ be a compact K\"ahler manifold and let $\rho:
\pi_1 (X)
\to S L_n (\Cal O_S)$.  Then the image under $\hat \rho_\ast$ of
$H_{2i -1} (X, \bbz)$ in $H_{2i -1} (B SL_n (\Cal O_S), \bbz)$ is torsion.
\endproclaim

Indeed, the values of the universal class $\hat c_i$ on the image of
$H_{2i-1} (X, \bbz)$ would
then lie in $(\bbc / \bbz)_{tors} = \bbq / \bbz$.

There were two fundamental ideas, which led us to the proof of this
result in the
case $i =2$.  The first was to make use of the Borel's theorem [4]
which states that the
stable real cohomology $H^\ast (BSL (\Cal O_S), \bbr)$ is multiplicatively
generated by $\sigma^\ast_k Vol_{2j-1}$.
Here $\{ \sigma_k \}$ is the maximal set of nonconjugate embeddings of
$F$ in $\bbc$ (we always
take $F$ big enough so that there would not be any real embeddings).  The
element $Vol_{2 j -1} \in H^{2 j -1} (B S L^\delta (\bbc), \bbr)$, the $\bbr$-
part of $\frac 1{\sqrt{-1}}\hat c_i$, is
the famous Borel regulator which we called the hyperbolic volume
regulator for $j = 2$.Recall that the $\bbr/\bbz$- part of $\hat c_i$ is the
(universal) Cheeger-Chern-Simons class.

Secondly, we developed a general approach to regulators of flat bundles
and gave a very precise geometrical recipe for their computation.  This
point of view results from our study in [21] and is relevant to the
approach of Corlette, c.f. [8],[9].

The fundamental analytic tool which was used in [22], was the
existence and rigidity for twisted
harmonic maps, a theory developed by Donaldson, Corlette, Sampson,
Carlson and Toledo.  The key
point was the use of the Sampson's degeneration theorem, which relies
heavily on the dimension
restriction $i =2$.  So to prove the Main Theorem 1.1 in it s full
generality, we need to invent some new ideas.
\pmf
{\bf 1.4} \ The strategy chosen in the present paper will be to use the
refined version of the Sampson's Theorem [25] and the deformation
theory of flat bundles over projective base, developed by Hitchin, Sampson,
Carlson -- Toledo, Goldman -- Millson and Simpson, for proving the following
key
result.  Two different proofs will be given.
\proclaim{Theorem}  Let $X$ be a compact K\"ahler manifold and let $\rho$ be as
in 1.1.  Then for any
$j \ge 2$ the
volume regulator $\hat \rho^\ast (Vol_{2j -1})
\in H^{2 j -1} (X, \bbr)$ is zero.
\endproclaim

For $j =2$ this is proved in [22], 4.3-4.5.  In the course of
the proof we will rely on the geometrical description
of the regulators, given in [22], 3.1-3.3.

I take the opportunity to thank all people, whose influence, personal or
through their work, determined my attitude to geometry.  I am deeply thankful
to Michael Gromov,
M. S. Narasimhan, Nigel Hitchin and Shicheng Wang.  I am particularly
grateful to Prof. Narasimhan for generous
sharing his insight on Donaldson's invariants, which became a starting
point in
my thinking on secondary invariants and their relation to moduli spaces
and algebraic $K$-theory.  Also, the lectures of Prof. Hitchin on Higgs
bundles during the Edinburgh Symposium in
1991 became a real clue for the approach presented here.I am equally grateful
to Carlos Simpson,Max Karoubi,H\'el\`ene Esnault  for interesting discussions
on regulators and flat bundles, and to the Editor, Pierre Deligne, for many
useful suggestions and improvements. Finally, I would
like to thank Professor Ilya Rips for the most attentive and
critical attitude to my work.  Special thanks are due to Charlotte Herman
for the excellent typing of the manuscript.This paper is dedicated to the
memory of my teacher, G.I.Katz.
\heading
2. Proof of the Main Theorem.
\endheading

Let $X, \rho$ be as in 1.1.  The proof of the theorem 1.1, and also
theorems 1.3 and 1.4 will be divided into the following steps.
\pmf
{\bf 2.1}  Consider the representation variety $V^{SL_n}_\Gamma$, where
$\Gamma = \pi_1 (X)$, and find a $\bar \bbq$-point $\bar \rho$ in the
component, containing $\rho$ (c.f. [22],2.2).  Then
$Ch _{2 i -1} (\bar \rho) = Ch _{2 i -1} ( \rho)$ ([7], Proposition 3.8 and
[22],5.16.1).
Moreover, we may assume that $\bar \rho$ is defined over some $\Cal O_S$
in a number field $F$ ( [22],2.2
).  Relabel $\bar \rho$ again by $\rho$.
\pmf
{\bf 2.2.} Consider a continuous map $\hat \rho : X \to B S L^\delta_n (\Cal
O_S)$, induced by
$\rho$.  There exists a universal regulator $\hat c_{2 i -1}
(B S L ({\Cal O}_S), \bbc / \bbz)$ such
that $\hat \rho^\ast \hat c_{2 i -1} = Ch _{2 i -1} ( \rho)$
( [1],[10],[12],[15]).
We wish to prove that $\hat \rho_\ast (H_{2 i -1} (X, \bbz)) \subseteq
H_{2 i -1} (BS L ({\Cal  O}_S), \bbz)$ is torsion.  Assuming this, we see
that the value of $\hat c_{2 i -1} \circ \hat \rho_\ast$ on any
class in $H_{2 i -1} (X, \bbz)$ will
lie in $(\bbc / \bbz)_{tors} = \bbq / \bbz$ as desired.

Now, $H_\ast (BS L (\Cal O_S), \bbz)$ is of finite type [5] and to
check the claim of Theorem 1.3 it suffices to test on
$\hat \rho_\ast (H_{2 i -1} (X, \bbz))$ the
generators of $H^{2i -1} (BSL (\Cal O_S), \bbr)$.  The theorem of
Borel [4] provides
a set of multiplicative generators for $H^\ast (BSL (\Cal O_S), \bbr)$
of the form
$\sigma^\ast_k Vol_{2 j -1}$, where $Vol_{2 j -1} \in H^{2 j -1}
(BSL (\bbc), \bbr)$ is the Borel
regulator and $\sigma_k : F \to \bbc$ form a maximal set of
nonconjugate embeddings.  We take
$F$ big enough so that there will be no real embeddings.  We claim
that for all
$k, j, \sigma^\ast_k Vol_{2j -1}$ vanishes on $\hat \rho_\ast
(H_{2 j -1} (X, \bbz))$, or,
equivalently, $Vol_{2j -1}$ vanishes an $\widehat {\sigma_k \circ \rho}_\ast
(H_{2 j -1} (X, \bbz))$ or further,
Vol$_{2j -1} (\sigma_k \circ \rho) = 0$ in $H^{2 j -1} (X, \bbr)$.
We refer to [22], for the detailed explanation of volume
invariants of representations.
\pmf
{\bf 2.3.}  We now see that the Theorem 1.1 and the Theorem 1.3. follows
from the Theorem 1.4.  applied to all Galois twists $\sigma_k \circ \rho_\ast$.
I call the reader's attention to the fact, that in order to prove the
 Theorem 1.1 even for an unitary representation (for which 1.4. trivially
holds) , we need to show the vanishing of volume invariant of
 non-unitary representations.

So we start to prove the Theorem 1.4. saying that Vol$_{2j -1} (\rho) = 0$.

Denote by ${\Cal F}$ the associated flat $SL_n ({\bbc}) / S U (n)$ - bundle
over $X$.  Recall that by [22],3.3, $Vol_{2j -1} (p)$ can be described
as follows.  Consider the canonical invariant $(2j -1)$ - form $\omega_{2 j
-1}$ on
$SL_n (\bbc) / S U (n)$.  Identifying the tangent space to the base point
$SU (n)$ with the space $p$ of traceless Hermitian matrices, this can be
written
as $\omega_{2j -1} (A_1, \dots A_{2 j -1}) = Alt (T r A_1 \dots A_{2 j -1} )$.
Next,
$\omega_{2 j -1}$ lifts to a closed $(2j-1)$ form on ${\Cal F}$ and using any
section, $s$, of ${\Cal F}$
we pull this lifted form down to $X$ to obtain a closed form, whose class,
denoted Bor$(\rho, \omega)$ in [22],3.2.2, gives $Vol_{2 j -1} (\rho) \in H^{2
j -1} (X, \bbr)$.
Since $SL_n (\bbc) / S U (n)$ is contractible, all sections are homotopic so
that the resulted class is independent on the choice of a section.

As in [22],4.5, we now take $s$ to be harmonic.  This is
possible by Donaldson - Corlette, at least if $\rho$ is semisimple.  Locally
we can view $s$ as a harmonic map $s: X \to S L_n (\bbc) / S U (n)$, defined
up to a left shift.  We claim that the pull down form $s^\ast \omega_{2 j -1}$
is identically zero.
Indeed, this will follow from the following result.
\pmf
{\bf 2.4 Theorem} (Sampson [25]).  Let $x \in X$.  Identify (uniquely up to a
conjugation) $T_{s(x)} {\Cal F}_x$ with $p$, where ${\Cal F}_x$ is a fiber
over $x$.  Then the image of $T^{1, 0}_x X$ under $D s \otimes \bbc: T^\bbc_x X
\to p^\bbc$ is
an abelian subspace of $p^\bbc$, and analogously for $T^{0, 1}_x X$.

\demo{Proof} \ For usual (untwisted) harmonic maps this is contained in [25].
The proof for harmonic sections is completely identical.
\pmf
{\bf 2.5} \ Now decomposing $s^\ast \omega_{2 j -1}$ to $(p, q)$-components we
see that
$s^\ast \omega_{2 j -1} = \sum^{2 j -1}_{i =0} c_i \nu_\cdot \lambda_i$, where
$c_i \in \bbc$ and for $v_1, \dots v_i \in T^{1, 0}_x X, \nu_i
(v_1 \dots v_i ) = Alt (Tr_\bbc \Big ( \Big ( D_x s \otimes \bbc (v_1)) \cdots
(D_x s \otimes \bbc (v_i) \Big )\Big )$, an similarly for $\lambda_i$.  Now 4.4
implies
immediately that $s^\ast \omega_{2 j -1} = 0$ for $j \ge 2$.In other words, if
$v_i=v_i'+v_i''$, such that $v_i'$ comutes, as well as $v_i''$, then $Alt
(Tr(v_1\cdots v_{2j-1})=0$.
\pmf
{\bf 2.6 Corollary} (cf [6]). Any continuous map of a
compact K\"ahler manifold to $\Gamma \setminus SL_n (\bbc) / S U (n)$ induces
trivial map in $H^i$ for $i \le m (SL_n (\bbc))$ (the Matsusima constant).
\demo{Proof} $H^\ast (SL_n (\bbc) / S U (n))$ is multiplicatively generated by
$\omega_{2 j -1}$ up to this range, by Matsusima ( c.f. [4]).
\pmf
{\bf 2.7} \ To deal with non-semisimple representations, consider the abelian
category $Mod(\bbc [ \pi_1 (X)])$ of finite-dimensional modules over $\pi_1
(X)$.  The
volume invariant $Vol_{2 j -1}$ is additive on exact sequences, in other words,
it defines a homomorphism $Vol_{2 j -1} : K_0 (Mod) \to H^{2 j -1} (X, \bbr)$.
So if
it vanishes an simple modules, it vanishes everywhere.
\pmf

\heading
3. Rudiments on Higgs Bundles and the Deformation Theory.
\endheading
\pmf
{\bf 3.1} Let $X$ be a compact K\"ahler manifold.  Consider a representation
$\rho : \pi_1 (X) \to SL_n (\bbc)$ and the corresponding flat vector
bundle $E_ \rho$.
If $\rho$ is unitary, equivalently, if there exists a parallel Hermitian
metric in $E_\rho$, then the behavior of the de Rham complex
$$ 0 \to \underline{\bbc} \otimes E_\rho \to C^\infty (E_\rho )
\to
\Omega^1 ( E_\rho) \to \cdots $$
is completely analogous to that in the untwisted case $\rho =1$.  Indeed, one
defines the Laplacian $\Delta : \Omega^i ( E_\rho ) \to \Omega^i
( E_\rho)$ using the flat connection
and the whole Hodge theory applies unchanged.  (In particular, for any $i$
odd, the dimension of $H^i ( \underline{\bbc} \otimes E_\rho)$ is
even and its positivity for some odd $i$ implies
$H^2 (\pi_1 (X) , \bbq) \neq 0)$.

For arbitrary representation one asks whether there exists a
distinguished Hermitian metric on $E_\rho$ such
that the scene for the Hodge theory is (partially) set.  The positive
answer for semisimple representation was given in the remarkable
series of papers by Nigel Hitchin [16],
Simon Donalson [11], Kevin Corlette [8], [9], and Carlos Simpson [27].  The
metric one chooses is just a harmonic section of the associated $SL_n (\bbc)/
S U (n)$ bundle.  Its existence is proved by the classical nonlinear heat
equation method of Eells-Sampson and the compactness of the
target is replaced by the sufficient twisting of the bundle, more
precisely, by the condition that the monodromy  action at
$S_\infty (SL_n (\bbc) / S U (n) )$ does not have fixed points, so
that it is energetically unprofitable for a minimizing
sequence of sections to leave all compact sets.
\pmf
{\bf 3.2.}  Recall that any harmonic map $f$ of a Riemannian surface $X$ to a
Riemannian manifold  $(M, g)$ gives rise to a holomorphic object - the
famous holomorphic
quadratic differential $(f^\ast g )^{2,0} \in H^0 (K^2)$. Completely
analogous,
a harmonic section of the $SL_n (\bbc) / S U (n)$ - bundle above gives
rise to a
holomorphic object - a holomorphic section $\theta$ of
$T^\ast_X \otimes End E_\rho$ for some new holomorphic structure in $E_\rho$.
The latter had got a name ``a Higgs field''
in the
N. J. Hitchin's ground - breaking paper [13].  Developing the gracious
theory of Narasimhan - Seshadri [19], Donaldson [11] and
Uhlenbeck  - Yau [29], which had established the bundle equivalence
flat unitary
$\leftrightarrow$ stable with vanishing $c_1, c_2$, Hitchin was able
to prove that any $\theta$-stable
holomorphic bundle over (one-dimensional) $X$ carries a ``Higgs metric'' $K$
for which
$F_K = - [ \theta, \theta^\ast]$.  Then an elementary computation shows that
the perturbed connection $\nabla_K + \theta + \theta^\ast$ is
flat, which makes the diagram
$$ \text{flat \ simple} \ \leftrightarrow \ \text{Higgs \ stable \ of \ degree
\ zero}, $$
a two-way highway.  Hitchin also discovered
a lot of additional geometrical structure in the moduli space of flat bundles
$\Cal G = Hom (\Gamma, SL_2 (\bbc))/ SL_2 (\bbc)$ (or rather its nonsingular
twisted counterpart), where
$\Gamma = \pi_1 (X)$ is the surface group.  In fact, the moduli space of stable
bundles $\Cal G = Hom (\Gamma, SU_2) / S U_2$ sits as a real symplectic
submanifold in $\Cal G$ with respect to the Goldman's symplectic structure.
Moreover,
$\Cal G$ is hyperk\"ahler and with respect to one special K\"ahler structure,
(which does not agree with the Goldman's structure)
induced by the conformal structure of $X, \Cal H \subset \Cal G$ is a complex
subvariety and $\Cal G$ contains a copy of $T^\ast_\bbc (\Cal H_{reg})$.
The leaves of $T^\ast_\bbc \Cal H$ over a smooth point (= a stable bundle $E$
over $X$) consists of all the Higgs
fields in $E$.  The quadratic map
$$ \det: T^\ast_\bbc (\Cal H) \to H^0 (K^2) $$
determines a bunch of holomorphic metrics on $\Cal H$ which Poisson
commute to each other.  The Lagrangian subvariety $\det^{-1} (pt)$ is
identified as the Primian of that double
covering of $X$, which makes the bundle $E$ to be an extension of
linear bundles (eigen bundles of $\theta$).
\pmf
{\bf 3.3.}  There is a clear and obvious need to extend this beautiful picture
to cover the representations of all K\"ahler groups.  In the part of the
equivalence relation this was realized by Simpson.  The
methods are largely the same, as far as the answer, which says
$$ \text{flat simple} \leftrightarrow \text{Higgs stable with} \ c_1, c_2 = 0.
\tag*$$

One remarkable feature of the key correspondence (*) is that the left side is,
so to say, much more subject to variations than the
right side.  For instance, the field $t \cdot \theta, t \in \bbc^\ast$ is again
a Higgs field and we recover following Simpson [27] a $\bbc^\ast$-action in
(*).  The
fixed points of this action correspond precisely to variations of complex Hodge
structure and the
Higgs operators $\theta (Z), Z \in T_x X$, are all nilpotent in this case.
Moreover,
Simpson establishes the following absolutely magnificent theorem.

\proclaim{Theorem} (Simpson).  Let $X$ be smooth projective.  Then any
representation $\rho:
\pi_1 (X) \to GL_n (\bbc)$ can be smoothly deformed to a variation
of a complex Hodge structure.
\endproclaim

Simpson gives a list of simple real Lie groups which may occur as the
(real) Zarishi closures of the monodromy subgroup $\rho (\pi_1 (X))
\subset GL_n (\bbc)$.  These groups, called the groups of Hodge type, are
divided to two subcategories:  the
symmetry groups of compact symmetric Hermitian manifolds and some other
groups,
including two infinite series and a number of exceptional groups.  By  a
result of Borel [3], any group of Hermitian type is really a Zariski
closure of a uniform lattice.  On
the other hand, I don't know wether all the groups of Hodge type which are
not of Hermitian type, actually occur as Zariski closures of the
monodromy representations.
\par
\newpage

\heading{4. Groups of Hodge type and the Cohomology of 
the Corresponding
Compact Symmetric Spaces.}
\endheading
\pmf

{\bf 4.1.}  This is the list of all simple noncompact groups of
Hodge type, by Simpson:
\par
{\settabs 3 \columns
\+$SU (p, q)$ &$Sp (n, \bbr)$\cr
\+$SO^\ast (2n)$ &$E_{6 ( - 14)}$\cr
\+$SO (p, 2)$ &$E_{7 ( - 25)}$\cr
\+ (all of Hermitian type)\cr
\+$SO (p, 2q), q \ge 2$ &$E_{8 (8)}$ \cr
\+$Sp (p,q)$ &$E_{8 (- 24)}$ \cr
\+$E_{6 (2)}$ &$F_{4 (4)}$\cr
\+$E_{7 (7)}$ &$F_{4 (- 20)}$ \cr
\+$E_{7 ( -5)}$ &$G_{2 (2)}$ \cr}
\pmf
{\bf 4.2.}  For any group $G$ from the list of 3.1. put $K \subset G$ to be
a maximal compact subgroup.  Then rank $K = \text{rank} \ G$ [27]
and $M = G / K$ is
a symmetric space.  The following table gives the Poincar\'e polynomials of the
dual compact
symmetric space $\hat M$ (c.f. Fomenko\break [14]).
\par
{\settabs 5 \columns
\+Group G &Subgroup $K$ &$G / K = M$ &$\hat M$ &$P_{\tilde M} (t)$ \cr
\+$SU(p. q)$ &$S U (p) \times S U (q)$ &$\frac{S U (p, q)}{S (U (p) \times
U (q))}$
&$\frac{S U (p + q)}{S (U (p) \times U (q))}$
&$\frac{(1 - t^{2 (q +1)} \dots (1 - t^{2(q+p)})}{(1 - t^2)
\dots (1 - t^{2p})}$\cr
\+$SO^\ast (2n)$ &$U (n)$ &$\frac{SO^\ast (2n)}{U (n)}$
&$\frac{SO (2n)}{U (n)}$ &$(1 + t^2) \dots (1 + t^{2n-2})$\cr
\+$SO (p, 2), p$ odd &$SO (2) \times SO (p)$
&$\frac{SO (p, 2)}{S O (2) \times SO (p)}$
&$\frac{SO (p +2)}{SO (2) \times SO (p)}$
&$\Big.  1 +t^2+ \dots + t^{2 p - 4}
 $ \cr
\+$SO(p, 2), p$ even &$SO(2) \times SO(p)$
&$\frac{SO(p, 2)}{SO(2)\times SO(p)}$
&$\frac{SO(p+2)}{SO(2) \times SO(p)}$
&$\Big (1 +t^2+\dots +t^{p - 2})\times$ \cr
\+ & & & &$\times (1 + t^{p - 2})$\cr
\+$Sp (n, \bbr)$ &$S U (n)$ &$\frac{Sp (n, \bbr)}{S U (n)}$
&$\frac{Sp (n)}{SU (n)}$
&$(1 + t^2) \dots (1 + t^{2n})$ \cr
\+$E_{6 ( -14)}$ &$SO (10) \times SO (2)$
&$\frac{E_{6 ( -14)}}{SO (10) \times
S O (2)}$ &$\frac{E_{6 ( - 78)}}{SO (10) \times SO (2)}$
&$(1 + t^2+ \dots + t^{16})\times$\cr
\+ & & & &$\times (1 + t^8 + t^{16})$\cr
\+$E_{7 ( - 25)}$ &$E_6 \times S O (2)$
&$\frac{E_{7 ( - 25)}}{E_6 \times SO (2)}$
&$\frac{E_{7 ( - 133)}}{E_6 \times S O (2)}$
&$(1 + t^2 + \dots + t^{26})\times$ \cr
\+ & & & &$\times (1 + t^{10})(1 + t^{18})$\cr
\+$SO (p, 2q), p$ odd &$SO (p) \times SO (2q)$
&$\frac{SO (p, 2q)}{SO (p) \times SO (2q)}$
&$\frac{SO (p + 2q)}{S O (p) \times SO (2 q)}$
&$\frac{( 1 - t^{2 (p + 1)} \dots (1 - t^{2 (p + 2q -1)})}{(1 - t^4)
\dots (1 - t^{4 ( q -1)})(1 - t^{2q})},$ \cr
\+$SO (p, 2q), p$ even &$SO (p) \times SO (2 q)$
&$\frac{SO (p, 2q)}{(SO (p) \times
SO (2q)}$ &$\frac{SO (p + 2 q)}{S O (p) \times S O ( 2 q)}$
&$\frac{(1 - t^{2p}) \cdots (1 - t^{2 (p + 2q -2)})}
{(1 - t^4) \cdots (1 - t^{4 (q - 1)})}$ \cr
\+$Sp (p, q)$ &$Sp (p) \times Sp (q)$
&$\frac{Sp (p, q)}{Sp (p) \times Sp (q)}$
&$\frac{Sp (p + q)}{Sp (p) \times Sp (q)}$
&$\frac{(1 - t^{4 (p + 1)}) \dots (1 - t^{4 (p + q)} )}{(1 - t^4)
\dots (1 - t^{4p})}$\cr
\+$E_{6 (2)}$ &$S U (6) \times SU (2)$ &$\frac{E_{6 (2)}}{S U (6) \times
SU (2)}$ &$\frac{E_{6 ( - 78)}}{S U (6) \times S U (2)}$
&$(1 + t^4 + \dots + t^{20})\times $\cr
\+ & & & &$\times (1 + t^6 + t^{12})(1 + t^8)$\cr
\+$E_{7 (7)}$ &$SU (8)$ &$\frac{E_{7 (7)}}{SU (8)}$ &$\frac{E_{7 ( -133)}}{S U
(8)}$
&$(1 + t^6 + \dots + t^{30})\times$ \cr
\+ & & & &$\times (1 +t^8+ t^{16}) \times$\cr
\+ & & & &$\times (1 + t^{10})(1 + t^{14})$\cr
\+$E_{7 ( - 5)}$ &$SO (12) \times S U (2)$
&$\frac{E_{7 ( - 5)}}{S O (12) \times S U (2)}$
&$\frac{E_{7 ( - 133)}}{S O (12) \times S U (2)}$
&$(1 + t^4 + \dots + t^{24})\times$\cr
\+& & & &$\times (1 + t^8 + t^{16}) \times $\cr
\+ & & & &$\times (1 + t^{12} + t^{24})$\cr
\+$E_{8 (8)}$ &$SO (16)$ &$\frac{E_{8 (8)}}{SO (16)}$
&$\frac{E_{8 ( -248)}}{S O (16)}$ &$(1 + t^8 + \dots +
t^{32})\times $\cr
\+& & & &$\times (1 + t^{12} + t^{24}) \times$ \cr
\+ & & & &$\times (1 + t^{16} + t^{32}) \times$ \cr
\+ & & & &$\times (1 + t^{20} + t^{40})$\cr
\+$E_{8 (-24)}$ &$E_7 \times S U (2)$ &$\frac{E_{8 ( - 24)}}{E_7 \times S U
(2)}$
&$\frac{E_{8 ( - 248)}}{E_7 \times S U (2)}$
&$(1 + t^4 + \dots + t^{36}) \times$\cr
\+ & & & &$\times (1 + t^{12} + t^{24} + t^{36})\times $\cr
\+ & & & &$\times (1 + t^{20} + t^{40})$ \cr
\+ $F_{4 (4)}$ &$Sp (3) \times S U (2)$ &$\frac{F_{4 (4)}}{Sp (3)
\times S U (2)}$
&$\frac{F_4}{Sp (3) \times S U (2)}$ &$(1 + t^4 + \dots + t^{20})\times$\cr
\+ & & & &$\times (1 + t^8)$ \cr
\+ $F_{4 (-20)}$ &$Spin (9)$ &$\frac{F_{4 ( -20)}}{Spin (9)}$
&$\frac{F_4}{Spin (9)}$ &$1 + t^8 + t^{16}$ \cr
\+$G_{2 (2)}$ &$SO (4)$ &$\frac{G_{2 (2)}}{SO (4)}$
&$\frac{G_2}{SO (4)}$ &$1 + t^4 + t^8$\cr}

We deduce the following proposition.
\pmf
\proclaim{4.3 Proposition} For any simple noncompact group $G$ of Hodge type,
the odd dimensional cohomology $H^{\text{odd}} ( \hat M) = 0$.  Consequently,
any $G$ - invariant form in $G / K$ is of even dimension.
\endproclaim

\demo{Proof}  The first statement follows from the inspection of 3.2. The
second is equivalent to it by duality (see Borel [4]).

\demo{4.4. Remark} The reader may wish to ask is there is a natural explanation
for the phenomenon, recovered in 3.3.  We refer to our paper [23] for the
relevant computation.
{\bf 4.5} \  We will give now a different proof of Theorem 1.4.
Using the Simpson's fundamental result we deform $\rho$ to a variation of
Hodge structure, which we call $\mu$.   By [7],
Vol$_{2 j -1} (\mu) = Vol_{2 j -1} (\rho)$.  So it is enough to show that
Vol$_{2j -1} (\mu) = 0$.
Let $G \subset SL_n (\bbc)$ be the Zariski closure of $\mu (\pi_1 (X))$, then
$G = \tilde K \times G_1 \times \dots \times G_n$, where $\tilde K$ is compact
and $G_i$
are from the list of 3.1.  Choose maximal compact subgroups $K_i \subset
G_i$ and $K \subset SL_n (\bbc)$ such that $K \cap G_i = K_i$, and $\tilde
K \subset K$.  Then we have a totally geodesic embedding
$$ G / K \subset SL_n (\bbc) / S U (n) $$

Now, let $\omega \in \Omega^{2 j -1} (SL_n (\bbc), SU (n))$ be an invariant
form such
that in terms of [22], Vol$_{2 j -1} (\mu) =
Bor (\omega, \mu)$.
Restricting an $G / K$ we get by naturality Vol$_{2 j -1} (\mu) = Bor
(\omega|_{G/ K}, \mu)$.
However, the form $\omega |_{G/ K}$ is an invariant form of odd dimension on
$G / K = \prod G_i / K_i$, which should be zero by 3.3.   Hence Vol$_{2 j -1}
(\mu) = 0$, which
proves theorem 1.4 and hence also theorem 1.3 and the Theorem 1.1.
\heading
5.  Representations in low Dimensions
\endheading
\pmf
{\bf 5.1.}  The aim of this section is to give a spirit of what kind
of behavior may
be expected from representations of K\"ahler groups in
$SL_n (\bbc)$, for $n$ small $(n \le 3)$.

By the fundamental result of Simpson, any representation may be
deformed to a variation of Hodge structure, so one should
begin to study these first.

Let $X$ be a K\"ahler manifold and let $\rho : \pi_1 (X) \to SL_2 (\bbc)$ is
a semisimple variation of Hodge structure .  The Zariski closure of
$\rho (\pi_1 (X))$ in $S L_2 (\bbc)$ is either a subgroup of $S U (2)$ or
$S L_2 (\bbr)$ or $P S U (1,1)
\approx S L_2 (\bbr)$.We have nothing to say in the first case.Consider
therefore a  representation in
$S U (1, 1)$, which
we relabel $\rho$.

Consider the flat bundle $E_\rho$.  By Simpson's theory, the
corresponding Higgs bundle decomposes as $L \oplus
M$ and the Higgs field $\theta$ lies in $\Omega^1 (X, Hom (L, M))$.
Observe that
$M = L^{-1}$ since there exists a parallel symplectic form in
$E_\rho$.  So $\theta \in \Omega^1 (X, M^2)$.  Next,
since $E_\rho$ is flat we get $0 = c_2 (E_\rho) = - c^2_1 (M)$.

So we get the following data:

(i) a line bundle $M$ with $c^2_1 (M) = 0$ and deg $M < 0$

(ii) a form $\theta \in \Omega^1 (X, M^2)$.

Conversely, given such data we construct a representation at
$\pi_1 (X)$ in $S L_2 (\bbr)$ by the Simpson's theory.  Observe that a
one-form with values in a line bundle gives rise to a codimension one
holomorphic foliation of $X$ (with singularities). In fact, Simpson and
Corlette proved that any nonrigid representation decomposes through an orbifold
curve.
\pmf
\demo{5.2 Example}  Suppose $X$ is a surface with
$H^{1,1} (X) \cap H^2 (X, \bbq)$
one-dimensional [24].  This is a generic condition.
Then by (i) we have $c_1 (M) = 0$ since the intersection form is positively
defined in the (one-dimensional) Hodge cycles subspace.  But this is
impossible, since deg $M$ should be negative.  So such surfaces do not admit
nontrivial homomorphisms $\rho : \pi_1 (X) \to S L_2 (\bbr)$.
\demo{5.3 Example}  Suppose $X$ is a Hilbert modular  surface, i.e.
$X = {\Cal H}^2 \times {\Cal H}^2 / \Gamma$,
where $\Gamma$ is a uniform irreducible lattice in $S L_2 (\bbr) \times S L_2
(\bbr)$.
Obviously $TX$ splits as $TX = L_1 \oplus L_2$.  Consider the line
bundle $L^{1/2}_i, i = 1, 2$.  It exists and may be
looked at as a spinor bundle, corresponding to a flat $\Cal H^2$-bundle over
$X$,
associated to defining representations $\rho_{1,2} : \Gamma \to S L_2 (\bbr)$.
(cf. [11]).  As such it carries the canonical Hermitian metric with nonpositive
curvature
form.  It is convenient to consider the local $\Cal H^2$-system over $X$,
associated
to $\rho_i, i = 1,2$, along with the canonical (developing) section.  Then
$L^{1/2}_i$ is
just $K^{-1/2}$, where $K$ is the canonical bundle of the fiber, looked at as a
Riemannian surface.  In particular, $(c_1 (L^{1/2}_i))^2 = 0$ and deg
$L^{1/2}_i < 0$.
Since rank $(S L_2 (\bbr) \times S L_2 (\bbr)) = 2, \Gamma$ is arithmetic, so
$X$ is projective.  The form $\theta_i \in \Omega^1 (X, L_i)$ is just the
projection $L_1 \oplus L_2 \to L_i$.  The two representations $\rho_i : \Gamma
\to S L_2 (\bbr)$ are rigid by Margulis.
\pmf
{\bf 5.4.}  Now we turn to the secondary class of $\rho$. This is just
$\lambda ((c^{Chow}_1 (M))^2)$, where $\lambda:Ch^i_0(X)\to H^{2i}_\cald
(X,\bbz(i))$ is the Deligne cycle class map and is torsion by [22], Main
Theorem, and  by 1.1.  If we remove
the condition (ii), we may ask the following question:
\pmf
\demo{Question}

Suppose $D$ is a divisor on $X$ such that $\deg D < 0$ and $[D]^2 = 0$ in
singular cohomology.  Consider $([D]^{Chow})^2 \in Ch^2_0 (X)$.  Is
$\lambda ([ D]^{Chow})^2$ always torsion?

We see that the answer is yes, if $T^\ast (X) \otimes \Cal L^{\otimes 2} (D)$
admits
a holomorphic section.
\pmf
{\bf 5.5.}  Now suppose $\Gamma$ is an arithmetic lattice in $S U (1,2)$ and $X
= B^2 / \Gamma$ is a hyperbolic
surface.  Let $\rho: \Gamma \to S U (1,2) \to S L (3, \bbc)$ be the defining
representation and let $E_\rho$ be the associated rank 3 flat vector
bundle.

We wish first to understand the underlying holomorphic structure of
$E_\rho$ (induced by the flat connection).  Realize $B^2$ as the quotient
of the cone $Q = ( | x_1 |^2 - | x_2|^2 - |x_3|^2 > 0)$ by the action of
$\bbc^\ast$.  We
will have then the canonically defined $SU (2, 1)$ - invariant tautological
line
bundle $L$ over $Q / \bbc^\ast = B^2$ (the analogue of ${\Cal O} (-1)$ over
$\bbp^2$) and an $S U (2, 1)$ - equivariant embedding $y : L \to B^n \times
\bbc^{n +1}$.  Descending to $X = B^2 / \Gamma$, we get a ``hyperplane''
bundle,
denoted also by $L$, and the embedding $y : L \to E_\rho$.  Moreover, the
quotient
$\bbc^{n +1} / \vp (L)$ is $T (B^n) \otimes L$, so we get an exact sequence
$$ 0 \to L \to E_\rho \to L \otimes T X \to 0 $$.
The extension class $\vp \in H^1 (Hom (L \otimes T X, L)) = H^1 (\Omega^1)$ is
proportional
to the K\"ahler class of $X$.  The corresponding Higgs bundle is a direct
sum $L \oplus L \otimes T X$ and $\theta$ is given by a matrix $\pmatrix
0 &0 \\
1 & 0 \endpmatrix$.    Next, $c^{Chow}_2 (E_\rho)$ is just $c^{Chow}_2 (X) -
\frac13
(c^{Chow}_1 (X))^2$.  So we get a rational class $B_\rho = \lambda
( c^{Chow}_2 (X) - \frac13 ( c^{Chow}_1 (X))^2) \in H^3 (X, \bbq / \bbz)$.
Observe
that $H^3 (X, \bbq / \bbz)$ may be nontrivial in many known cases, but it seems
difficult to compute $B_\rho$.
\pmf
{\bf 5.6} \  We will use the description of 5.5 to give a proof of a version
of the following famous theorem of Yau.

\proclaim{Theorem} Let $X$ be a compact K\"ahler surface such that
{\text i)} \ $\Omega^2 (X)$ is positive and
{\text ii)} \ $c^2_1 (X) = 3 c_2 (X)$.

Then $X$ is a ball quotient.
\endproclaim

We will replace (i) by a weaker assumption that
deg$(\Omega^2 (X)) = - c_1 (X) \omega > 0$ but
demand $T X$ to be stable (e.g. there are no holomorphic foliations on $X$).We
valso make a technical assumption that $\frac{1}{3}c_1(X)\in H^2(X,\bbz)$.
Form a Higgs bundle $L \oplus L \otimes T X$ where
$L^{\otimes 3} = \Omega^2 (X)$, and
$\theta$ is $\pmatrix 0 &0 \\
1 &0 \endpmatrix$.  We claim it is stable.  In fact, any $\theta$ - invariant
subbundle is a subbundle of $L \otimes T X$, say $L \otimes E$ where
$E \subset T X$.  Compute deg$(L \otimes E) = \deg L$ rank
$E + \deg E = - \frac13 \deg T X\cdot$ rank
$E + \deg E$, so $\deg (L \otimes E) < 0$ is equivalent to
$\frac{\deg E}{rank E} < \frac13 \deg T X$ which
is weaker then stability of $TX$.
Since $c_1 (L \oplus L \otimes TX) = 0$ and
$c_2 ( L \oplus L \otimes TX) = 0$, we come to a variation of
Hodge structure, whose monodromy should be a Hodge type
subgroup of $SL_3 (\bbc)$.  It is easy to see that there is the only
opportunity of $SU (2, 1)$.  So we get a representation $\rho : \pi_1 (X) \to S
U (2, 1)$ and
$E_\rho$ corresponds to $L \oplus L \otimes TX$.

Next, find a harmonic metric in $E_\rho$.  This is a harmonic section of
the associated $S U (2, 1) / S U (2) = B^2$ - bundle over $X$, say, $s$.  By
Sampson-Corlette - Carlson - Toledo, we have that either $Ker(Ds|_x)$ is a
nonzero
complex subspace of $T_x X$, or $s$ is holomorphic.  But the field $\theta$ is
just $(D s \otimes \bbc)^{1,0}$, and viewed as a map $\theta_x : T_x X^{1, 0}
\to End (L \oplus L \otimes TX)$ is of
maximal rank everywhere.  This rules out the first possibility, so $s$ is
holomorphic, and the same
argument shows that rank$_\bbc D s =2$ everywhere, so we can pull the metric of
$B^2$ down to $X$ using $s$ and hence $X$ has a hyperbolic metric.

This type of argument was first used in [21] to prove the Goldman's theorem.
Observe that for surfaces of general type the
condition (i) is 5.6 can be removed by Bogomolov-Miyaoka.
\pmf
{\bf 5.7} \ We wish to find a similar characterization of Hilbert modular
surfaces.
Recall that any K\"ahler manifold diffeomorphic to such surface is actually
bi-holomorphic to it.
We suggest the following purely geometrical description.
\proclaim{Theorem} \ Let $X$ be a compact K\"ahler surface such that
{\text i)} \ The signature $\sigma (X) = 0$ and $\chi (X) \neq 0$
{\text ii)} \ There exist two transversal holomorphic foliations on $X$ with
nonpositive normal bundles.

Then $X$ is biholomorphic to a Hilbert modular surface.
\endproclaim

\demo{Proof} \ We have $T X = L_1 \oplus L_2$ and the both $L_i$ are
nonpositive.  Compute $0 = p_1 (T X) = p_1 (L_1) + p_1 (L_2) = c^2_1
(L_1) + c^2_1 (L_2)$, so $c^2 (L_1) = c^2_1 (L_2) = 0$.
Moreover, $\chi (X) = c_2 (X) = c_1 (L_1) c_1 (L_2) \neq 0$, so $c_1 ( L_i)
\neq 0$.Since $c_1(L_i)$ lies in the Hodge cycle subspace of $H^{1,1}$, we have
by the  Hodge index theorem deg $L_i = c_1 (L_i) \omega \neq 0$
and since $L_i$ is nonpositive, deg $L_i < 0$.  We will make a technical
assumption $c_1 (L_i)$ is even, i.e.
$\frac12 c_1 (L_i) \in H^2 (X, \bbz)$.  Then $L_i = M^2_i$ for some line bundle
$M_i$ and by 5.1  we get two representations $\rho_{1, 2}: \pi_1 (X) \to SL_2
(\bbr)$.  Let $s_i$ be the
harmonic section of the corresponding that
${\Cal H}^2$ - bundle, then $L_1 = Ker (D s_2)$ and $L_2 = Ker (D s_1)$ so
$(s_1, s_2)$ is of maximal rank everywhere.  This imposes the locally
homogeneous structure on $X$, as above.
\pmf
{\bf 5.8} \ We close with applications of Higgs theory to higher Milnor
inequalities,
refering to [21] for a general survey.  Let first $X$ be a compact Reimann
surface of
genus $g > 1$ and let $\rho : \pi_1 (X) \to S U (n, 1)$ be a representation.
The
classifying space $B S U^\delta (\infty, 1)$ carries a special element $\hat c
\in H^2 (B S U^\delta (\infty, 1))$, introduced and studied by
Morita [18], Toledo [28], Corlette [9] and the author [21].  The
pull-back $\hat p^\ast c \in H^2 (X, \bbr)$ is a secondary class of $X$, and
the number $(\hat p^\ast c, [X])$ is
called the degree of $\rho$, denoted deg$\rho$.  Then we have the following
beautiful result.
\proclaim{Theorem} (Milnor [17], Hitchin [16] $(n = 1)$, Toledo [28]).  For
any $X, \rho$, we have $| \deg \rho | \le g - 1$.
\endproclaim

\demo{Proof} As in [19], we may suppose $\rho$ to be irreducible.  Let ${\Cal
F}$
be the associated flat $B^n$-bundle over $X$ and let $s$ be
a harmonic section.  Recall that there is a natural $SU (n, 1)$-equinvariant
line bundle $L$ over $B^n$, equipped with the canonical connection.  Let $N$ be
the
pull-back bundle over $X$ by $s$ with the holomorphic structure, induced
by the pull-back connection.  Let $M$ be the pull-back of the tangent bundle to
$B^n$, again with the Levi-Civita connection and the
induced holomorphic
structure (over $X$).  The flat bundle $E_\rho$, as a complex bundle, is
an extension of $M \otimes N$ by $N$, whereas the corresponding Higgs bundle is
$N \oplus M \oplus N$ and the
$\theta$-field is given by $\pmatrix 0 &\ast \\
(D s \otimes \bbc)^{1, 0} &0 \endpmatrix$.  Moreover, the degree deg$\rho$ is
just deg$N$.

Since $X$ is a curve, the image of $N$ by $\theta$ extends to a holomorphic
subbundle $P$ of $M \otimes N$, so that $0 \neq \theta \in H^0 (\Omega^1
\otimes Hom (N, P))$ (the case $\theta = 0$
is trivial).  Hence deg$P \ge \deg N - ( 2 g - 2)$.  The subbundle $N \oplus P$
is
$\theta$-invariant, so by the Higgs stability, deg$P + \deg N < 0$ which means
$\deg N < g - 1$, if $n > 1$.  Changing the
orientation on $X$, we see that also $- \deg \rho < g - 1$, so $| \deg \rho | <
g - 1$.  The
equality, for $\rho$ irreducible is possible only if $n =1$ and $\Omega^1
\otimes Hom (N, P)$ trivial, which gives the standard
Higgs bundle $K^{+ 1/2} \oplus K^{- 1/2}, \theta = \pmatrix 0 &a \\
1 &0 \endpmatrix$, corresponding to the embedding $\pi_1 (X) \to S U (1, 2) \to
S U (n, 1)$, by Hitchin [16].
\pmf
{\bf 5.9} \ The case of higher-dimensional $X$ is conceptually simpler because
of
Siu's superrigidity, which forces $s$ to be holomorphic.  We just state the
corresponding
inequalities, leaving the proofs to the reader.

\proclaim{Theorem} \ Let $X = B^m / \Gamma$ be a ball quotient and let $\rho :
\Gamma \to S U (n, 1)$ be a
representation.  Then $\Big | \Big ( \hat p^\ast c \cdot \omega^{m -1}, [ X]
\Big ) \Big | \le \Big ( \omega^m, [ X] \Big ) = Vol X$.
\endproclaim

Moreover [9], $\Big | \Big ( \hat p^\ast c)^m, [X] \Big | \le Vol X$.

It is remarkable that one can prove the parallel statement in real
hyperbolic geometry, namely, that for any representation $\rho : \pi_1 (X) \to
PSO (n, 1)$ of the
fundamental group of a compact hyperbolic manifold, one gets $Vol (\rho) \le
Vol (X)$.  We refer to
[22] for the details.
\pbf
\centerline{References}
\item{1.} A. A. Beilinson, {\it Higher regulators and values of $L$-functions},
J. Soviet Math.
{\bf 30} (1985), 2036--2070.
\item{2.} S. Bloch, {\it Applications of the dilogarithm function in algebraic
$K$-theory and
algebraic geometry}, Int. Symp. on Alg. Geom., Kyoto, 1977, 103--114.
\item{3.} A. Borel, {\it Compact Clifford-Klein forms of symmetric spaces},
Topology,
{\bf 2} (1963), 111-122.
\item{4.} A. Borel, {\it Stable real cohomology of arithmetic groups}, Annals
Sci. Ec. Norm. Sup\'er,
{\bf 7} (1974), 235--272; II, Prog. Math. Boston, {\bf 14} (1981), 21--55.
\item{5.} A. Borel, J. - P. Serre, {\it Cohomogie d'immeubles et de groupes
$S$-arithm\'etique}, Topology {\bf 15} (1976), 211--232.
\item{6.} J. Carlson, D. Toledo, {\it Harmonic mappings of K\"ahler manifolds
to locally symmetric spaces}, Publ. Math. IHES, {\bf 69} (1989), 173--201.
\item{7.} S.-S. Chern, Simons, {\it Characteristic  forms and geometric
invariants},
Annals of Math., {\bf 99} (1974), 48--69.
\item{8.} K. Corlette, {\it Flat $G$-bundles with canonical metrics},
J. Diff. Geom., {\bf 28} (1988), 361--382.
\item{9.} K. Corlette, {\it Rigid representations of K\"ahlerian fundamental
groups}, J. Diff. Geom., {\bf 33} (1991), 239--252.
\item{10. } P.Deligne, {\it Th\'eorie de Hodge, I,II,III}, Actes Congr. Intern.
Math. Nice, 1970,425--430; Publ.Math.IHES {\bf 40} (1971), 5--58,{\bf
44}(1974),5-77.
\item{11.} S. K. Donaldson, {\it Twisted harmonic maps and self-duality
equations},
Proc. London Math. Soc., {\bf 55} (1987), 127--131.
\item{12. }H.Esnault, E.Vieweg, {\it Deligne--Beilinson cohomology}, in:
Beilinsom's Conjectures on Special values of L-functions,M.Rappoport,
N.Schappacher, P.Schneider, ed.,Acad.Press,1988,43--91.
\item{13. } H.Esnault, V.Srinivas, {\it Chern classes of vector bundles with
holomorphic connections on a complete smooth complex variety},J.Diff.Geom, {\bf
36 }, (1992), 257--267.  .
\item{14.} A. T. Fomenko, {\it Variational Methods in Topology},
Nauka, 1982 (Russian).
\item{15. }H.Gillet, C.Soul\'e, {\it Arithmetic Chow groups and differential
characters,}, in: ``Algebraic $K$-theory; connections with Geometry and
Topology'',Kluwer, 1989, 29--68.
\item{16.} N. J. Hitchin, {\it The self-duality equations on a
Riemann surface},
Proc. London Math. Soc., {\bf 55} (1987), 59--126.
\item{17.} J. Milnor, {\it On the existence of a connection with curvature
zero}, Comm.Math.Helv., {\bf 32}(1958), 215--223.
\item{18.} S. Morita, {\it Characteristic classes of surface bundles} Inv.
Math. {\bf 100} (1985), 3--16.
\item{19.} M. S. Narasimhan, C.S. Seshadri,  {\it  Stable and unitary bundles
on a compact Riemann surface}, Ann. of Math. {\bf 82} (1965), 540--564.
\item{20.} D. Ramakrishnan, {\it Regulators, algebraic cycles and values
of $L$-function},
Contemp. Math. {\bf 85} (1989), 183--310.
\item{21.} A.Reznikov,  {\it Harmonic maps, hyperbolic cohomology and higher
Milnor inequalities}, Topology {\bf 32},(1993),899--907.
\item{22.} A. Reznikov, {\it Rationality of secondary classes}, preprint
(April, 1993).
\item{23.} A. Reznikov, {\it Morse theory and combinatorial identities},
preprint (1992).
\item{24.} J. Rogawski, to appear.
\item{25.} J. Sampson, {\it Applications of harmonic maps to K\"ahler
geometry},
Contemp. Math. {\bf 49} (1986), 125--133.
\item{26.} C. Soul\'e, {\it Connexions et classes caract\'eristiques de
Beilinson}, Contemp. Math., {\bf 83}
(1989), 349--376.
\item{27.} C. Simpson, {\it Higgs bundles and local systems}, Publ. Math.
IHES, {\bf 75} (1992), 5--95.
\item{28.} D. Toledo, {\it Representations of surface group in complex
hyperbolic space}, J. Diff. Geom. {\bf 29} (1989), 125--133.
\item{29.} K. Uhlenbeck, S.-T. Yau, {\it On the existence of
Hermitian-Yang-Mills connections in
stable vector bundles}, Comm. Pure and Appl. Math., {\bf 39} (1986), 257--293.
\pbf
\item{} Institute of Mathematics
\item{} The Hebrew University
\item{} Givat Ram, 91904 Jerusalem
\item{} ISRAEL
\item{} email:simplex@sunrise.huji.ac.il
\bye